\documentclass[a4paper,11pt]{article}
\usepackage{latexsym}
\newtheorem{thm}{Theorem}[section]
\newtheorem{lem}[thm]{Lemma}

\newtheorem{cor}[thm]{Corollary}
\newtheorem{prop}[thm]{Proposition}

\newtheorem{rem}[thm]{Remark}

\title{A new inequality for the Hermite constants}
\author{Roland Bacher\footnote{Support 
from the Swiss National Science Foundation is
gratefully acknowledged.}}
\date{}

\begin{document}
\maketitle


\par {\it Abstract: We describe continuous
increasing functions $C_n(x)$ such that 
$\gamma_n\geq C_n(\gamma_{n-1})$ where $\gamma_m$ is Hermite's
constant in dimension $m$. This inequality yields a 
new proof of the Minkowski-Hlawka bound 
$\Delta_n\geq \zeta(n)2^{1-n}$ for the maximal density 
$\Delta_n$ of $n-$dimensional lattice-packings.}
\footnote{Math. class.: 10E05, 10E20. Keywords: Lattice packing,
Hermite constant, inequality}
 
\section{Introduction and main results}
                          
We denote by $\hbox{min}(\Lambda)=\hbox{min}_{\lambda\in\Lambda
\setminus\{0\}}\langle \lambda,\lambda\rangle$ the minimum 
(defined as the squared Euclidean 
length of a shortest non-zero element) of 
an $n-$dimensional lattice $\Lambda\subset {\mathbf E}^n$ in
the Euclidean vector-space ${\mathbf E}^n$ 
and define the {\it density} of $\Lambda$ by
$$\Delta(\Lambda)=\sqrt{\frac{(\hbox{min}\ \Lambda)^n}{4^n\ \hbox{det}
\ \Lambda}}V_n$$
where $V_n=\frac{\pi^{n/2}}{(n/2)!}$ stands through 
the whole paper for the volume of the
$n-$dimensional unit-ball in ${\mathbf E}^n$. The density
$\Delta(\Lambda)$ is the proportion 
of vo\-lume occupied by a maximal 
open Euclidean ball embedded in the flat torus ${\mathbf
  E}^n/\Lambda$ with volume $\sqrt{\det(\Lambda)}$ and having
a shortest closed geodesic of length $\sqrt{\hbox{min}(\Lambda)}$. 
The largest density $\Delta_n=\Delta(\Lambda_n)$
achieved by an $n-$dimensional lattice $\Lambda_n$ is called 
the {\it maximal density} in dimension $n$. 
Related constants are the {\it maximal center density} 
$\delta_n=\Delta_n/V_n$ and the 
{\it Hermite constant} $\gamma_n=4\delta_n^{2/n}$  in 
dimension $n$. The sequence $\gamma_1,\gamma_2,\dots$ of
Hermite constants satisfies for $n\geq 3$ {\it Mordell's inequality}
$$\gamma_n\leq \gamma_{n-1}^{(n-1)/(n-2)}$$
which yields an upper bound for $\gamma_n$ (if $n\geq 3$)
in terms of $\gamma_{n-1}$. Our main result is a
complementary inequality bounding $\gamma_n$ from below in terms
of $\gamma_{n-1}$. For the convenience of the reader we state it in
three equivalent ways, either in terms of densities $\Delta_m$,
center-densities $\delta_m$ or Hermite constants $\gamma_m$ in
dimension $m$. It involves the {\it M\"obius function}
$\mu:{\mathbf N}_{>0}\longrightarrow{\mathbf Z}$ defined by
$\mu(l)=(-1)^a$ for a natural integer $l\in {\mathbf N}$ which 
is a product of $a$ distinct primes and by $\mu(l)=0$ 
if $l$ is divisible by the square of a prime number. 

\begin{thm} \label{Hermiteineq} (i) The maximal densities
$\Delta_{n-1}$ and
  $\Delta_{n}$ of lattice-packings in dimensions $n-1$ and $n\geq 2$ 
satisfy the inequality
$$2^{n-1}\Delta_{n-1}
\sum_{k=1}^{\lfloor2\Delta_nV_{n-1}/(\Delta_{n-1}V_n)\rfloor}
\left(\sum_{l\vert k}\frac{\mu(l)}{l^{n-1}}\right)\left(
1-\left(
\frac{k\Delta_{n-1}V_n}{2\Delta_nV_{n-1}}\right)^2\right)^{(n-1)/2}\geq 1$$
where the sum $\sum_{l\vert k}$ is over all positive integral
divisors $l\in{\mathbf N}$ of the natural integer $k$. 

\ \ (ii) The maximal center densities $\delta_{n-1}$ and
  $\delta_{n}$ of lattice-packings in dimensions $n-1$ and $n\geq 2$ 
satisfy the inequality
$$2^{n-1}\delta_{n-1}\frac{\pi^{(n-1)/2}}{((n-1)/2)!}
\sum_{k=1}^{\lfloor2\delta_n/\delta_{n-1}\rfloor}
\left(\sum_{l\vert k}\frac{\mu(l)}{l^{n-1}}\right)\left(
1-\left(\frac{k\delta_{n-1}}{2\delta_n}\right)^2\right)^{(n-1)/2}\geq 1.$$

\ \ (iii) The Hermite constants $\gamma_{n-1}$ and $\gamma_n$ in
dimensions $n-1$ and $n\geq 2$ satisfy the inequality
$$\frac{\pi^{(n-1)/2}}{((n-1)/2)!}\sum_{k=1}^{\lfloor
\gamma_n^{n/2}/\gamma_{n-1}^{(n-1)/2}\rfloor}
\left(\sum_{l\vert k}\frac{\mu(l)}{l^{n-1}}\right)\left(
\gamma_{n-1}
-k^2\left(\frac{\gamma_{n-1}}{\gamma_n}\right)^n\right)^{(n-1)/2}\geq 1.$$
\end{thm}

Consider the function
$$F_n(x,y)=\sum_{k=1}^{\lfloor \sqrt x\  y\rfloor}\sum_{l\vert k}
\frac{\mu(l)}{l^{n-1}}\left(x-\left(\frac{k}{y}\right)^2\right)^{(n-1)/2},
$$
defined for $x,y>0$. It is increasing in $y$, strictly increasing 
in $y$ for $y>\frac{1}{\sqrt{x}}$, continuous and extends 
continuously to $y=0$ since $F_n(x,y)=0$ for $y\leq\frac{1}{\sqrt x}$. 
There exists thus a continous function
$x\longmapsto Y_n(x)$ such that 
$$F_n(x,Y_n(x))=\frac{((n-1)/2)!}{\pi^{(n-1)/2}}$$
for all $x>0$.
Set $$C_n(x)=\hbox{sup}_{0<\xi\leq x}(\xi\ (Y_n(\xi))^{2/n}).$$

\begin{thm} \label{mainineq}
Let $\tilde \gamma_{n-1}\leq \gamma_{n-1}$ be a lower
  bound for Hermite's constant in dimension $n-1$. Then
$$\gamma_n\geq C_n(\tilde \gamma_{n-1})\ .$$
\end{thm} 

Analogous results hold of course for
$\Delta_n=2^{-n}\gamma_n^{n/2}V_n$ 
and $\delta_n=2^{-n}\gamma^{n/2}$ related to Hermite's constant 
$\gamma_n=4(\Delta_n/V_n)^{2/n}=4\delta_n^{2/n}$.

\begin{rem} (i) The inequality of Theorem \ref{mainineq}
is tight for $n=2$. For $n=3$, we get from 
$\delta_2=1/2\sqrt{3}$ the lower bound $0.1695\leq \delta_3=1/4\sqrt
2\sim
0.1768$. For $n=9$, the known value
$\delta_8=1/16$ gives the lower bound $\delta_9\geq 0.0388$
(a lattice with center-density $0.0442$ is known),
for $n=25$ the known value $\delta_{24}=1$ coming from the Leech
lattice (see Cohn-Kumar, \cite{CK1} and \cite{CK2}) 
yields $\delta_{25}\geq 0.657$ (a lattice with center-density
$0.707$ is known). 

\ \ (ii) The above examples show that our inequality is better than
the tri\-vial inequality 
$$\delta_n\geq \sqrt{\frac{\mu_{n-1}^n}{4^n
\det(\Lambda_{n-1}\oplus \mu_{n-1}{\mathbf Z})}}=
\frac{\delta_{n-1}}{2}$$
obtained by considering the orthogonal sum $\Lambda_{n-1}\oplus \mu_{n-1}
{\mathbf Z}$ of
a densest $(n-1)-$dimensional lattice $\Lambda_{n-1}$ with
minimum $\mu_{n-1}=\hbox{min}_{\lambda\in\Lambda_{n-1}
\setminus\{0\}}\langle \lambda,\lambda\rangle$.

\ \ (iii) The factor 
$$\sum_{l\vert k}\frac{\mu(l)}{l^{n-1}}=\prod_{p\hbox{ prime}, p\vert
  k}
\left(1-\frac{1}{p^{n-1}}\right)$$
yields only a minor improvement for huge $n$ and is the analogue of a
standard trick leading to the factor $\zeta(n)$
in the Minkowski-Hlawka bound $\Delta_n\geq \frac{\zeta(n)}{2^{n-1}}$.

\ \ (iv) The application $x\longmapsto C_n(x)$ of Theorem 
\ref{mainineq} is strictly increasing for $x$ huge enough.
The computations of Section \ref{secMink} show in fact
that one has $C_n(x)=x(Y_n(x))^{2/n}$ (at least for huge $n$)
for all values of $x$ which are of interest.
\end{rem}

\begin{rem} \label{Marin}
Starting with the inequality of assertion (ii) in Theorem 
\ref{Hermiteineq}, and using $\sum_{l\vert k}\frac{\mu(l)}{l^{n-1}}=
\prod_{p\vert k}\left(1-\frac{1}{p^{n-1}}\right)<1$ 
(where the product is over all prime divisors of $k$), 
A. Marin pointed out the easy inequalities
$$\begin{array}{l}
\displaystyle 1\leq
2^{n-1}\delta_{n-1}\frac{\pi^{(n-1)/2}}{((n-1)/2)!}
\sum_{k=1}^{\lfloor 2\delta_n/\delta_{n-1}\rfloor}
\left(\sum_{l\vert k}\frac{\mu(l)}{l^{n-1}}\right)\left(
1-\left(\frac{k\delta_{n-1}}{2\delta_n}\right)^2\right)^{(n-1)/2}\\
\displaystyle\quad \leq  
2^{n}\delta_{n}V_{n-1}\sum_{k=1}^{\lfloor 2\delta_n/
\delta_{n-1}\rfloor}\left(
1-\left(\frac{k\delta_{n-1}}{2\delta_n}\right)^2\right)^{(n-1)/2}
\frac{\delta_{n-1}}{2\delta_n}\\
\displaystyle\quad \leq 2^n\delta_n V_{n-1}\int_{0}^1
\left(1-x^2\right)^{(n-1)/2}dx=2^n\delta_n\frac{V_n}{2}=2^{n-1}\Delta_n
\end{array}$$
which show $\Delta_n\geq \frac{1}{2^{n-1}}$. This is, up to the factor
$\zeta(n)$, the Minkowski-Hlawka bound for $\Delta_n$.
\end{rem}

The following, technically more involved result yields a slightly
better lower bound:

\begin{thm} \label{Minkowskibound} For all $\epsilon>0$, there 
exists $N$ such that 
$$\Delta_n\geq\frac{1-\epsilon}{2^n\ \sum_{k=1}^\infty e^{-k^2\pi}}\sim
(1-\epsilon)\ 23.1388\ 2^{-n}$$
for all $n\geq N$.
\end{thm}

\begin{rem}
Theorem \ref{Minkowskibound} is slightly better that
the Minkowski-Hlawka bound which shows the
existence of lattices
with density at least $\zeta(n)\ 2^{1-n}$, cf. formula (14) in  \cite{CS},
Chapter 1. The best known bound for densities achieved by
lattice packings
(together with a very nice proof) seems to be due to Keith Ball
and asserts the existence of $n-$dimensional lattices with density 
at least $2(n-1)2^{-n}\zeta(n)$, see \cite{Ball}.
Previously, similar results where proven by Rogers and
Davenport-Rogers, see \cite{R} and \cite{DR}. Somewhat related are 
also \cite{KLV} and \cite{GZ} which describe elementary constructions
of dense sphere and lattice packings.
\end{rem}

The paper is organized as follows:

Section \ref{sectiondef} introduces notations and summarizes
for the convenience of the reader a few well-known facts
on lattices. It contains also an easy (and seemingly not
very well-known) result on integral sublattices which are orthogonal
to a non-zero integral vector in ${\mathbf Z}^{n+1}$.

In Section \ref{sectionmu} we define of $\mu-$sequences
which are the main tool of this paper. Theorem \ref{majorationI} 
of this Section gives a quantitative (and somewhat technical)
statement for extending a suitable finite $\mu-$sequence 
$(s_0,\dots,s_{n-1})$ to a $\mu-$sequence
$(s_0,\dots,s_{n-1},s_n)$. The associated
$n-$dimensional lattice $(s_0,\dots,s_n)^\perp\cap
{\mathbf Z}^{n+1}$ is obtained by a close analogue of  lamination
with respect to its $(n-1)-$dimensional sublattice
$(s_0,\dots,s_{n-1})^\perp\cap{\mathbf Z}^{n}$.

Theorem \ref{majorationI}
is the central result of this paper since it implies easily
Theorem \ref{Hermiteineq} and Theorem \ref{mainineq}
as shown at the end of Section \ref{sectionmu}. 
The proof of Theorem \ref{Minkowskibound} is more technical 
and given in Section \ref{secMink}.

Section \ref{sectioneasy} states and proves a weaker 
and easier statement
than Theorem \ref{majorationI}. 
Although not necessary for the other parts of 
the paper, this section describes a fairly elementary and
almost effective method for constructing dense lattices. It contains
moreover the essence of the main idea for proving 
Theorem \ref{majorationI}.  

Section \ref{sectionproof} describes the proof of Theorem 
\ref{majorationI}.

Section \ref{secMink} is devoted to the proof of 
Theorem \ref{Minkowskibound}.

Section \ref{sectioncomments} contains a few final comments
and remarks.


\section{Definitions}\label{sectiondef}

All facts concerning lattices needed in the sequel are collected 
in this Section for the convenience of the reader,
see  \cite{CS} and \cite{M} for more on lattices and
lattice-packings. 

An {\it $n-$dimensional lattice} is a discret-cocompact 
subgroup $\Lambda$ 
of the $n-$dimensional Euclidean vector space ${\mathbf E}^n$.
Denoting by $\langle\ ,\ \rangle$ the scalar product and choosing
a ${\mathbf Z}-$basis $b_1,\dots,b_n$ of a lattice 
$\Lambda=\oplus_{j=1}^n {\mathbf Z} b_j$, 
the positive definite symmetric matrix $G\in {\mathbf R}^{n\times n}$
with coefficients
$$G_{i,j}=\langle b_i,b_j\rangle$$
is a {\it Gram matrix} of $\Lambda$. Its determinant $\det(G)$,
called the {\it determinant}
of $\Lambda$, is independent of the choosen basis $b_1,\dots b_n$
and equals the squared volume of the flat torus ${\mathbf E}^n/\Lambda$.
The {\it norm} of a lattice vector
$\lambda\in\Lambda$ is defined as $\langle \lambda,\lambda\rangle$
and equals thus the square of the Euclidean norm $\sqrt{\langle
\lambda,\lambda\rangle}$. A lattice
$\Lambda$ is {\it integral} if all scalar products
$\{\langle\lambda,\mu\rangle\ \vert\ \lambda,\mu\in\Lambda\}$
are integral. An integral lattice
of determinant $1$ is {\it unimodular}. An Euclidean 
lattice $\Lambda$ is unimodular if and only if every group homomorphism
$\varphi:\Lambda\longrightarrow{\mathbf Z}$ is of the form
$\varphi(v)=\langle v, w_\varphi\rangle$ for a suitable fixed 
element $w_\varphi\in\Lambda$. The {\it minimum} 
$$\hbox{min}\ \Lambda\ =\hbox{min}_{\lambda\in\Lambda\setminus\{0\}}
\langle \lambda,\lambda\rangle$$
of a lattice $\Lambda$ is the norm
of a shortest non-zero vector in $\Lambda$. The {\it density}
$\Delta(\Lambda)$ and the {\it center-density}
$\delta(\Lambda)$ of an $n-$dimensional lattice $\Lambda$ are defined as
$$\Delta(\Lambda)=\sqrt{\frac{(\hbox{min}\ \Lambda)^n}{4^n\ \hbox{det}
\ \Lambda}}\ V_n\qquad \hbox{and}\qquad 
\delta(\Lambda)=\sqrt{\frac{(\hbox{min}\ \Lambda)^n}{4^n\ \hbox{det}
\ \Lambda}}$$
where $V_n=\pi^{n/2}/(n/2)!$ denotes the volume of the $n-$dimensional
unit-ball in ${\mathbf E}^n$.
These two densities are proportional for a given fixed dimension $n$ 
and $\Delta(\Lambda)$
measures the (asymptotic) proportion of space occupied by the {\it sphere
packing of $\Lambda$} 
obtained by centering $n-$dimensional Euclidean balls
of radius $\sqrt{\hbox{min}\ \Lambda/4}$ at all points of $\Lambda$.

Given an $n-$dimensional lattice $\Lambda\subset {\mathbf E}^n$ the subset
$$\Lambda^\sharp=\{x\in{\mathbf E}^n\ \vert\ \langle x,\lambda\rangle\in{\mathbf Z}
\quad \forall \lambda\in\Lambda\}$$
is also a lattice called the {\it dual lattice} of $\Lambda$. 
The scalar product induces a natural bijection between 
$\Lambda^\sharp$ and the set of homomorphisms
$\Lambda\longrightarrow {\mathbf Z}$.
A lattice $\Lambda$ is integral if and only if 
$\Lambda\subset \Lambda^\sharp$. For an
integral lattice, the {\it determinant group} $\Lambda^\sharp/\Lambda$ is
a finite abelian group consisting of $(\hbox{det}\ \Lambda)$ elements.

A sublattice $M\subset \Lambda$ is {\it saturated}
if $\Lambda/M$ is without torsion or equivalently
if $M=(M\otimes_{\mathbf Z}{\mathbf R})\cap \Lambda$. A sublattice
$\Lambda\subset M$ is thus saturated if and only if 
$M$ is a direct factor of the additive group $M$.

The following result  is well-known (cf. Chapter I, 
Proposition 9.8 in \cite{M}):

\begin{prop} \label{proporthogreseaux} 
Let $M$ and $N$ be two saturated sublattices of dimension $m$ and $n$ 
in a common unimodular lattice $\Lambda$ of dimension $m+n$ such that
$M$ and $N$ are contained in orthogonal subspaces.

Then the two determinant groups
$M^\sharp/M\hbox{ and }N^\sharp/N$
are isomorphic. In particular, the
lattices $M$ and $N$ have equal determinants.
\end{prop}

\noindent{\bf Proof} Since $\Lambda=\Lambda^\sharp$ is unimodular,
orthogonal projection $\Lambda\longrightarrow M^\sharp$
yields a surjective homomorphism from $\Lambda$ onto $M^\sharp$
with kernel $(N\otimes_{\mathbf Z} {\mathbf R})\cap \Lambda$
coinciding with $N$ since $N$ is saturated. This shows
$M^\sharp\sim \Lambda/N$ and thus $M^\sharp/M\sim \Lambda/(M\oplus
N)$. Exchanging the role of $M$ and $N$ implies the result.
\hfill$\Box$

Two lattices $\Lambda$ and $M$ are {\it similar}, 
if there exists a bijection 
$\Lambda\longrightarrow M$ which extends to an Euclidean
similarity from $\Lambda\otimes_{\mathbf Z}{\mathbf R}$ 
to $M\otimes_{\mathbf Z}{\mathbf R}$. The set of 
similarity classes of lattices is endowed with a natural topology:
a neighbourhood of an $n-$dimensional lattice $\Lambda$ is given by all
lattices having a Gram matrix in ${\mathbf R}_{>0}\ V(G)$
where $V(G)\subset {\mathbf R}^{n\times n}$ 
is a neighbourhood of a fixed Gram matrix $G$ of $\Lambda$.

Similar lattices have identical densities and the density function
$\Lambda\longmapsto \Delta(\Lambda)$ is continuous with respect
to the natural topology on similarity classes.

Consider the set ${\cal L}_n=\{\Lambda(s)\ \vert\ s\in
{\mathbf N}^{n+1}\setminus\{0\}\}$ of all $n-$dimensional 
integral lattices of the form
$$\Lambda(s)=\{z\in{\mathbf Z}^{n+1}\ \vert\ \langle z,s\rangle=0\}$$
where $s\in{\mathbf N}^{n+1}\setminus\{0\}$ is a non-zero vector
of length $n+1$ with non-negative integral coordinates.

\begin{prop} \label{Lestdense} The set ${\cal L}_n$ is dense in the set of
similarity classes of $n-$dimensional Euclidean lattices.
\end{prop}

There are thus lattices in ${\cal L}_n$ with densities arbitrarily
close to the ma\-ximal density $\Delta_n$ of $n-$dimensional 
lattices. 

\noindent {\bf Proof of Proposition \ref{Lestdense}} 
Given a Gram matrix $G=\langle b_i,b_j\rangle$ 
of an $n-$di\-men\-sional lattice 
$\Lambda=\oplus_{j=1}^n{\mathbf Z}b_j$, Gram-Schmidt orthogonalization
of the ${\mathbf Z}-$basis $b_1,\dots,b_n$
(with respect to the Euclidean scalar product) 
yields a matrix factorization
$$G=L\ L^t$$
where $L=(l_{i,j})_{1\leq i,j\leq n}$ is an invertible 
lower triangular matrix.

Choose a large real number $\kappa>0$ and consider the integral lower 
triangular matrix
$\tilde L(\kappa)$ whose coefficients $\tilde l_{i,j}\in{\mathbf Z}$ 
satisfy
$$\vert \tilde l_{i,j}-\kappa l_{i,j}\vert\leq 1/2$$
and are obtained by rounding off each coefficient of $\kappa L$ 
to a nearest integer.

Define the integral matrix
$$B(\kappa)=\left(\begin{array}{ccccccc}
\tilde l_{1,1}&1&0&0&\dots\\
\tilde l_{2,1}&\tilde l_{2,2}&1&0&\\
\vdots&&\ddots&\ddots\\
\tilde l_{n,1}&\tilde l_{n,2}&\dots&\tilde l_{n,n}&1
\end{array}\right)$$
of size $n\times (n+1)$ with coefficients
$$b_{i,j}=\left\lbrace\begin{array}{ll}
\tilde l_{i,j}&\hbox{if }j\leq i\\
1&\hbox{if }j=i+1\\
0&\hbox{otherwise .}\end{array}\right.$$ 

The rows of $B(\kappa)$ 
span an integral
sublattice $\tilde \Lambda(\kappa)$ of dimension $n$ in ${\mathbf Z}^{n+1}$.
Moreover, the lattice $\tilde \Lambda(\kappa)$ is saturated since
deleting the first column of $B(\kappa)$ yields an integral
unimodular square matrix of size $n\times n$. 
The special form of $B(\kappa)$ shows that there exists an integral 
column-vector
$$v(\kappa)=\left(\begin{array}{c}
1\\-\tilde l_{1,1}\\
\tilde l_{1,1}\tilde l_{2,2}-\tilde l_{2,1}\\
\vdots\end{array}\right)\in {\mathbf Z}^{n+1}$$
such that $B(\kappa)v(\kappa)=0$.   
We have thus 
$$\tilde \Lambda(\kappa)=
v(\kappa)^\perp\cap{\mathbf Z}^{n+1}\subset {\mathbf E}^{n+1}\ .$$

Since $\hbox{lim}_{\kappa\rightarrow\infty}\frac{1}{\kappa}\ B(\kappa)$
is given by the matrix $L$ with an extra column of zeros appended,
we have
$$\hbox{lim}_{\kappa\rightarrow\infty}\frac{1}{\kappa^2}\ B(\kappa)
(B(\kappa))^t=G$$
and the lattice $\frac{1}{\kappa}\tilde \Lambda(\kappa)$
converges thus to the lattice $\Lambda$ for $\kappa\rightarrow \infty$.
Considering the integral vector $s=(s_0,s_1,\dots)\in{\mathbf
  N}^{n+1}$ defined by
$s_i=\vert v(\kappa)_{i+1}\vert$ for $i=0,\dots,n$, we get an integral
lattice 
$$\{z=(z_0,\dots,z_n)\in{\mathbf Z}^{n+1}\ \vert\ \langle 
z,s\rangle=0\}$$
of ${\mathcal L}_n$ 
which is isometric to $\tilde \Lambda(\kappa)$.\hfill$\Box$

\section{$\mu-$sequences}\label{sectionmu}

Let $\mu\geq 2$ be a strictly positive integer. 
A {\it $\mu-$sequence} is a finite or infinite
sequence $s_0=1,s_1,s_2,\dots$ of $(l+1)$ 
strictly positive integers such that 
the $n-$dimensional lattice
$$\Lambda_n=\{(z_0,z_1,\dots,z_n)\in{\mathbf Z}^{n+1}\ \vert \ 
\sum_{k=0}^n s_kz_k=0\}=(s_0,\dots,s_n)^\perp\cap{\mathbf Z}^{n+1}$$
has minimum $\geq \mu$ for all $n\geq 1$ which make sense 
(ie. for $n\leq l$ if the sequence 
$(s_0,s_1,\dots,s_{l})$ has finite length $l$). 
The letter $\mu$ in a $\mu-$sequence stands for minimum and should
not be confused with the Moebius-function, (unfortunately also)
denoted $\mu:{\mathbf N}\longrightarrow \{-1,0,1\}$.

Since $\Lambda_n$ is
saturated in ${\mathbf Z}^{n+1}$ 
by construction and orthogonal to the $1-$dimensional 
saturated lattice ${\mathbf Z}(s_0,\dots,s_n)\subset {\mathbf Z}^{n+1}$,
Proposition \ref{proporthogreseaux} shows that we have 
$\hbox{det}(\Lambda_n)=\sum_{k=0}^n s_k^2$. We get thus a 
lower bound for the density 
$$\Delta(\Lambda_n)=\sqrt{\frac{(\hbox{min}\ \Lambda_n)^n}{4^n\ \hbox{det}
\ \Lambda_n}}V_n\geq \sqrt{\frac{\mu^n}{4^n\ \sum_{k=0}^n s_k^2}}V_n$$
of the $n-$dimensional lattice $\Lambda_n$ associated to a
$\mu-$sequence $(s_0,\dots,s_n,\dots)$. This lower bound is an
equality except if the sequence 
$(s_0,\dots,s_{n})$ is a $(\mu+1)-$sequence.

\begin{rem} We hope that the double meaning of $\mu$ will not confuse
the reader: $\mu(l)\in\{-1,0,1\}$ denotes always the M\"obius function
of a natural integer $l$ while $\mu$ or $\mu_1,\mu_2,\dots$
stands for natural integers.
\end{rem}

\begin{rem} \label{rem1} (i) The condition $s_0=1$ ensures
that $(s_0,\dots,s_n)\subset{\mathbf Z}^{n+1}$ generates a 
saturated $1-$dimensional sublattice and will be useful for 
proving Lemma \ref{lemmaXk}. It can however be weakened by
requiring that $s_0,\dots,s_n$ are without common non-trivial
divisor. Lemma \ref{lemmaXk} (which applies to a sequence
of $\mu-$sequences) remains valid if $s_0$ is uniformly bounded.

It is of course 
also possible (but not very useful) to consider sequences
with coefficients in ${\mathbf Z}$.

\ \ (ii) Any subsequence $s_{i_0}=s_0,s_{i_1},s_{i_2},\dots$
of a $\mu-$sequence is again a $\mu-$sequence and permuting the
terms of a $\mu-$sequence by a permutation fixing $s_0$
yields of course again a $\mu-$sequence. 

\ \ (iii) Lattices associated to $\mu-$sequences are generally
neither perfect nor eutactic (cf. \cite{M} for definitions).
Their densities can thus generally be improved by suitable deformations.
\end{rem}

\begin{thm}\label{majorationI} Let $\mu_1,\mu_2,\dots$
be a strictly increasing sequence of natural integers $2\leq
\mu_1<\mu_2<\dots$. Suppose that we have finite
$\mu_r-$sequences $(s(\mu_r)_0,\dots,s(\mu_r)_{n-1})$ with existing
limit-density
$$\tilde
\Delta_{n-1}=\hbox{lim}_{r\rightarrow \infty}
\frac{\mu_r^{(n-1)/2}}{\sqrt{4^{n-1}\sum_{i=0}^{n-1}
    s(\mu_r)_i^2}}V_{n-1}>0$$
for the sequence of lattices 
$(s(\mu_r)_0,\dots,s(\mu_r)_{n-1})^\perp\subset {\mathbf Z}^n$.

Let $\sigma_n$ be a positive real number such that
$$2^{n-1}\tilde\Delta_{n-1}\sum_{k=1}^A
\left(\sum_{l\vert k}\frac{\mu(l)}{l^{n-1}}\right)
\left(1-k^2\left(2^{n-1}\tilde \Delta_{n-1}\frac{V_n}
{V_{n-1}}\sigma_n\right)^2\right)^{(n-1)/2}<1$$
where 
$$A=\left\lfloor \frac{2^{1-n} V_{n-1}}{V_n\tilde \Delta_{n-1}\sigma_n}
\right\rfloor.$$
Then there exists a natural integer $R$ such that for all $r\geq R$, 
the $\mu_r-$se\-quence 
$(s(\mu_r)_0,\dots,s(\mu_r)_{n-1})$
can be extended to a $\mu_r-$sequence
$(s(\mu_r)_0,\dots,s(\mu_r)_{n})$ satisfying 
$0<s(\mu_r)_n<\sigma_n\mu_r^{n/2}V_n$.
\end{thm}

The proof of Theorem \ref{majorationI} will be given in section
\ref{sectionproof}. Assuming Theorem \ref{majorationI}, 
we proceed now to prove Theorem \ref{Hermiteineq}.

\subsection{Proof of Theorem \ref{Hermiteineq} and \ref{mainineq}} 

Suppose that the inequality 
of assertion (i) does not hold for some natural integer $n$. By
Proposition \ref{Lestdense} we can find a sequence
of finite $\mu_r-$sequences $(s(\mu_r)_0,\dots,s(\mu_r)_{n-1})$
(with $\mu_r\rightarrow\infty$) such that 
$$\hbox{lim}_{r\rightarrow\infty}
\frac{\mu_r^{(n-1)/2}}{\sqrt{4^{n-1}\sum_{i=0}^{n-1}
    s(\mu_r)_i^2}}V_{n-1}=\Delta_{n-1}\ .$$
Consider the function
$$x\longmapsto 2^{n-1}\Delta_{n-1}\sum_{k=1}^{A(x)}
\left(\sum_{l\vert k}\frac{\mu(l)}{l^{n-1}}\right)\left(1-k^2\left(2^{n-1}
\Delta_{n-1}\frac{V_n}
{V_{n-1}}x\right)^2\right)^{(n-1)/2}$$
where $A(x)=\lfloor \frac{2^{1-n} V_{n-1}}{V_n\Delta_{n-1}x}
\rfloor$ for $x>0$. We claim that this function is continuous and
strictly decreasing in $x$. Indeed, for increasing $x\in
[0,\frac{2^{1-n}V_{n-1}}{kV_n\Delta_{n-1}}]$,
a summand
$$\left(\sum_{l\vert k}\frac{\mu(l)}{l^{n-1}}\right)\left(1-k^2\left(2^{n-1}
\Delta_{n-1}\frac{V_n}
{V_{n-1}}x\right)^2\right)^{(n-1)/2}$$ 
decreases continuously from 
$\sum_{l\vert k}\frac{\mu(l)}{l^{n-1}}=\prod_{p\vert k}
\left(1-\frac{1}{p^{n-1}}\right)>0$ to $0$. Such a summand
disappears if it becomes zero and their number 
$A(x)=\lfloor \frac{2^{1-n} V_{n-1}}{V_n\Delta_{n-1}x}\rfloor$
is finite and decreases with $x$.

This shows that we can choose a positive real number
$\sigma_n<\frac{1}{2^n\Delta_n}$ such that we have the inequalities
$$\begin{array}{l}
\displaystyle 2^{n-1}\Delta_{n-1}
\sum_{k=1}^{\lfloor2\Delta_nV_{n-1}/(\Delta_{n-1}V_n)\rfloor}
\left(\sum_{l\vert k}\frac{\mu(l)}{l^{n-1}}\right)
\left(1-\left(
\frac{k\Delta_{n-1}V_n}{2\Delta_nV_{n-1}}\right)^2\right)^{(n-1)/2}<\\
\displaystyle \quad < 
2^{n-1}\Delta_{n-1}\sum_{k=1}^A
\left(\sum_{l\vert k}\frac{\mu(l)}{l^{n-1}}\right)\left(1-k^2\left(2^{n-1}
\Delta_{n-1}\frac{V_n}
{V_{n-1}}\sigma_n\right)^2\right)^{(n-1)/2}<1\end{array}$$
where $A=\lfloor \frac{2^{1-n} V_{n-1}}{V_n\Delta_{n-1}\sigma_n}
\rfloor$.

Applying Theorem \ref{majorationI} and extracting a suitable 
subsequence from
$\mu_1,\mu_2,\dots$, we can suppose that all sequences
$(s(\mu_r)_0,\dots,s(\mu_r)_{n-1})$ can be extended to $\mu_r-$sequences
$(s(\mu_r)_0,\dots,s(\mu_r)_{n})$ with 
$$\hbox{lim}_{r\rightarrow \infty}\frac{s(\mu_r)_n}{\mu_r^{n/2}}=\alpha\leq
\sigma_nV_n\ .$$
Using 
$$\hbox{lim}_{r\rightarrow\infty}\frac{1}{\mu_r^{n-1}}\sum_{i=0}^{n-1}
s(\mu_r)_i^2=\frac{V_{n-1}^2}{4^{n-1}\Delta_{n-1}^2}$$
we have
$$\hbox{lim}_{r\rightarrow\infty}\frac{1}{\mu_r^n}\sum_{i=0}^n s(\mu_r)_i^2=
\hbox{lim}_{r\rightarrow\infty}\frac{1}{\mu_r}\ 
\frac{V_{n-1}^2}{4^{n-1}\Delta_{n-1}^2}+\alpha^2=\alpha^2$$
and get the existence of a sequence of $n-$dimensional lattices 
$$(s(\mu_r)_0,\dots,s(\mu_r)_n)^\perp\subset{\mathbf Z}^{n+1}$$
with limit-density
$$\hbox{lim}_{r\rightarrow\infty}
\sqrt{\frac{\mu_r^{n}}{4^n\sum_{i=0}^ns(\mu_r)_i^2}}V_n=
\frac{1}{2^n\alpha}V_n\geq\frac{1}{2^n\sigma_nV_n}V_n>\Delta_n$$
in contradiction with maximality of $\Delta_n$.\hfill$\Box$

{\bf Proof of Theorem \ref{mainineq}} Choose $\xi<x$ such that $C_n(x)=\xi\
Y_n(\xi)^{2/n}$. The inequality $\xi\leq \gamma_{n-1}$ implies the
existence of an $(n-1)$--dimensional lattice with density $
2^{-(n-1)}\xi^{(n-1)/2}V_{n-1}$.
Theorem \ref{majorationI} implies the existence of $n-$dimensional
lattice with density $2^{-n}(C_n(x))^{n/2}V_n$.
\hfill$\Box$

\section{An easy crude bound for the lexicographically
first $\mu-$sequence}
\label{sectioneasy}

\begin{thm} \label{mainA} Given an integer $\mu\geq 2$ there 
exists an increasing
$\mu-$sequence $s_0=1\leq s_1\leq \dots$ such that 
$$s_n\leq 1+\sqrt{\mu-2}\sqrt{\mu-1+n/4}^n\frac{\sqrt{\pi}^n}{(n/2)!}
\leq \sqrt{\mu}\sqrt{\mu+n/4}^n\frac{\sqrt{\pi}^n}{(n/2)!}$$
for all $n\geq 1$.
\end{thm}

The proof of Theorem \ref{mainA} is elementary and
consists of an analysis of the  \lq\lq greedy algorithm''
which constructs the first $\mu-$sequence with respect to
the lexicographic order on sequences. An easy analysis shows that the
lexicographically first $\mu-$sequence satisfies the first inequalities
of Theorem \ref{mainA}. The greedy algorithm, although very simple, 
is unfortunately useless for practical purposes.

A $\mu-$sequence satisfying the inequalities of Theorem \ref{mainA}
yields already rather dense lattices as shown by the next result.

\begin{cor} \label{corA} For any $\mu\geq 2$, there exists a 
$\mu-$sequence
$(s_0,s_1,\dots,s_n)\in{\mathbf Z}^{n+1}$ such that the density of
the associated lattice $\Lambda_n=(s_0,\dots,s_n)^\perp\cap{\mathbf Z}^{n+1}$
satisfies
$$\Delta(\Lambda_n)
\geq \frac{(1+n/(4\mu))^{-n/2}}{2^{n}\sqrt{(n+1)\mu}}\ .$$
\end{cor}

\begin{rem} \label{rem2} Taking $\mu\sim n^2/4$ 
we get for large $n$ the existence of $n-$dimensional lattices
with density $\Delta$ roughly at least equal to
$$\frac{1}{2^{n-1}\ n\ \sqrt{(n+1)\ e}}$$
which is reasonably close to the Minkowski-Hlawka bound
$\Delta_n\geq \zeta(n)\ 2^{1-n}$.
\end{rem}

\begin{lem} \label{lemZn} The standard Euclidean lattice ${\mathbf Z}^n$ contains at
most
$$2\sqrt{\mu+n/4}^n\frac{\pi^{n/2}}{(n/2)!}$$
vectors of (squared Euclidean) norm $\leq \mu$.
\end{lem}

{\bf Proof} We denote by $$B_{\leq \sqrt \rho}(x)=\{z\in{\mathbf E}^n\ 
\vert\ \langle z-x,z-x\rangle\leq \rho\}$$
the closed Euclidean ball with radius $\sqrt{\rho}\geq 0$
and center $x\in{\mathbf E}^n$. Given $\sqrt{\mu},\sqrt{\rho}\geq 0$
and $x\in B_{\leq \sqrt\mu}(0)$, the closed half-ball
$$\{z\in{\mathbf E}^n\ \vert\ \langle z,x\rangle\leq \langle x,x\rangle\}
\cap B_{\leq \sqrt\rho}(x)$$
(obtained by intersecting the closed affine halfspace $H_x=
\{z\in{\mathbf E}^n\ \vert\ \langle z,x\rangle\leq \langle x,x\rangle\}$ with the Euclidean ball $B_{\leq \sqrt \rho}(x)$ centered at $x\in
\partial H_x$) is contained in $B_{\leq \sqrt{\mu+\rho}}(0)$.

Since the regular standard cube
$$C=[-\frac{1}{2},\frac{1}{2}]^n\subset {\mathbf E}^n$$
of volume 1 is contained in a ball of radius $\sqrt{n/4}$ 
centered at the origin, 
the intersection
$$(z+C)\cap\{x\in {\mathbf E}^n\ \vert \ \langle x,x\rangle\leq \mu+n/4\}=
(z+C)\cap B_{\leq \sqrt{\mu+n/4}}(0)$$
is of volume at least $1/2$ for any element $z\in {\mathbf E}^n$ of
norm $\langle z,z\rangle \leq \mu$.

Since integral translates of $C$ tile ${\mathbf E}^n$, we have
$$\frac{1}{2}\sharp\{z\in{\mathbf Z}^n\ \vert\ \langle z,z\rangle\leq \mu\}\leq
\hbox{Vol}\ \{x\in{\mathbf E}^n\ \vert\ \langle x,x\rangle\leq \mu+n/4\}
\ .$$
Using the fact that the unit ball in Euclidean 
$n-$space has volume $\pi^{n/2}/(n/2)!$ (cf. 
Chapter 1, formula 17 in \cite{CS}) we get the result. \hfill $\Box$

{\bf Proof of Theorem \ref{mainA}} For $n=0$, the first inequality boils
down to $s_0=1\leq 1+\sqrt{\mu-2}$ and holds for
$\mu\geq 2$. Consider now for $n\geq 1$ a $\mu-$sequence
$(s_0,\dots,s_{n-1})\in {\mathbf N}^n$.

Introduce the set
$${\cal F}_n=\{(a,k)\in{\mathbf N}^2\ \vert\ \exists\  z=
(z_0,\dots,z_{n-1})\in
{\mathbf Z}^n\setminus\{0\}\hbox{ such that}$$
$$ak=\langle z,(s_0,\dots,s_{n-1})\rangle\hbox{ and }
\langle z,z\rangle+k^2<\mu\}\ .$$

Since $\Lambda_{n-1}$ has minimum $\geq \mu$, the equality
$\langle z,(s_0,\dots,s_{n-1})\rangle=0$ implies $\langle z,z\rangle
\geq \mu$ for $z\in {\mathbf Z}^n\setminus\{0\}$.
This shows that we have $a,k>0$ for $(a,k)\in{\mathcal F}_n$.

Since for a given pair of opposite non-zero 
vectors $\pm z\in{\mathbf Z}^n$ with norm
$0<\langle z,z\rangle<\mu$ there are at most 
$\sqrt{\mu-1-\langle z,z\rangle}\leq \sqrt{\mu-2}$ strictly 
positive integers
$k$ such that $\langle z,z\rangle+k^2<\mu$, such a pair $\pm z$
of vectors contributes at most $\sqrt{\mu-2}$
distinct elements to ${\cal F}_n$. The cardinality $f_n=
\sharp({\cal F}_n)$
of ${\cal F}_n$ is thus bounded by
$$f_n\leq \sqrt {\mu-2}\frac{\sharp\{z\in{\mathbf Z}^n\ \vert \ 
0<\langle z,z\rangle\leq \mu-1\}}{2}\leq 
\sqrt {\mu-2}\sqrt{\mu-1+n/4}^n\frac{\pi^{n/2}}{(n/2)!}$$
where the last inequality follows from Lemma \ref{lemZn}. 
There exists thus a strictly positive integer 
$$s_n\leq f_n+1\leq 1+\sqrt {\mu-2}\sqrt{\mu-1+n/4}^n\frac{\pi^{n/2}}{(n/2)!}$$
such that $(s_n,k)\not\in {\cal F}_n$ for all $k\in {\mathbf N}$. 
The strictly positive integer $s_n$
satisfies the first inequality of Theorem \ref{mainA} and it is 
straightforward to check that the $n-$dimensional lattice
$$\Lambda_n=\{z\in{\mathbf Z}^{n+1}\ \vert\ \sum_{i=0}^n s_iz_i=0\}$$
has minimum $\geq \mu$. This shows the first inequality. 
Choosing for $s_n$ the smallest strictly positive integer such that
$(s_n,k)\not\in {\mathcal F}_n$ for all $k\in{\mathbf N}$ and
iterating this construction yields clearly an increasing 
$\mu-$sequence.

The second inequality
$$1+\sqrt{\mu-2}\sqrt{\mu-1+n/4}^n\frac{\sqrt{\pi}^n}{(n/2)!}
\leq \sqrt{\mu}\sqrt{\mu+n/4}^n\frac{\sqrt{\pi}^n}{(n/2)!}$$
of Theorem \ref{mainA} boils down to 
$$1\leq \sqrt 2\sqrt{2+n/4}^n\ \frac{\sqrt{\pi}^n}{(n/2)!}$$
for $\mu=2$. This inequality is clearly true since the $n-$dimensional
Euclidean ball of radius $\sqrt{2+n/4}$ has volume
$\sqrt{2+n/4}^n\ \frac{\sqrt{\pi}^n}{(n/2)!}$ and contains the regular
cube $[-\frac{1}{2},\frac{1}{2}]^n$ of volume 1.

For $\mu\geq 3$ we have to establish the inequality
$\Phi(1)-\Phi(0)\geq 1$ where
$$\Phi(t)=\sqrt{\mu-2+2t}\sqrt{\mu-1+t+n/4}^n\frac{\sqrt{\pi}^n}{
(n/2)!}\ .$$
We get thus
$$\begin{array}{lcl}
\Phi(1)-\Phi(0)&\geq&\hbox{inf}_{\xi\in(0,1)}\Phi'(\xi)\\
&\geq&
\frac{1}{\sqrt\mu}\sqrt{\mu-1+n/4}^n\frac{\sqrt{\pi}^n}{(n/2)!}
+\frac{n}{2}\sqrt{\mu-2}\sqrt{\mu-1+n/4}^{n-2}\frac{\sqrt{\pi}^n}{(n/2)!}\ .
\end{array}$$

For $n=1$ and $\mu\geq 2$ we have 
$$\Phi(1)-\Phi(0)\geq\sqrt{1-\frac{1}{\mu}}\frac{\sqrt{\pi}}{\sqrt{\pi}/2}\geq 
\frac{2}{\sqrt 2}>1\ .$$

For $n\geq 2$ and $\mu\geq 3$ we get
$$\Phi(1)-\Phi(0)\geq \sqrt{2+n/4}^{n-2}\frac{\sqrt{\pi}^{n-2}}{((n-2)/2)!}\ 
\pi$$
and the right-hand side equals $\pi>1$ for $n=2$. 
For $n>2$, the right hand side equals $\pi$ times the volume of the
$(n-2)-$dimensional ball of radius $\sqrt{2+n/4}$ containing
the regular cube $[-\frac{1}{2},\frac{1}{2}]^{n-2}$ of volume $1$.
The second inequality follows.
\hfill $\Box$

{\bf Proof of Corollary \ref{corA}} Theorem \ref{mainA} shows the
existence of a $\mu-$sequence $(s_0=1,\dots,s_n)$ satisfying
$$s_0,\dots,s_n\leq \sqrt\mu\sqrt{\mu+n/4}^n\frac{{\sqrt \pi}^n}{(n/2)!}\ .$$
This shows for the lattice 
$\Lambda_n=(s_0,\dots,s_n)^\perp\cap{\mathbf Z}^{n+1}$ the inequality
$$\det\Lambda_n=
\sum_{i=0}^n s_i^2\leq (n+1)\mu(\mu+n/4)^n\frac{\pi^n}{((n/2)!)^2}=
(n+1)\mu(\mu+n/4)^n V_n^2$$
and implies
$$\Delta(\Lambda_n)\geq\sqrt{\frac{\mu^n}{4^n(n+1)
\mu(\mu+n/4)^n V_n^2}}V_n$$
which proves Corollary \ref{corA}.\hfill $\Box$

\section{Proof of Theorem \ref{majorationI}}\label{sectionproof}

The main idea for proving Theorem \ref{majorationI} is to get rid 
of a factor 
$\sqrt\mu$ when computing an upper bound $f$ for
the size of the finite set $\mathcal F$ considered in the
proof of Theorem \ref{mainA}. This is possible since the volume 
of the standard unit-ball of large dimension concentrates along 
linear hyperplanes. During the proof, we use for simplicity the slightly 
abusive notation $\mu=\mu_k$ and
$(s_0,\dots,s_{n})=(s(\mu_k)_0,\dots,s(\mu_k)_n)$. 
Since $\mu$ belongs to the strictly increasing integral sequence
$\mu_1<\mu_2<\dots$ tending to infinity, 
we consider sequences in the $\mu\rightarrow\infty$
limit. This allows us to neglect boundary effects when replacing
counting arguments by volume-computations.

In the sequel we write
$$g(x)\sim_{x\rightarrow\alpha}h(x)\ ,\hbox{ respectively }g(x)\leq_{x\rightarrow\alpha}h(x)\ ,$$
for $$\hbox{lim}_{x\rightarrow\alpha}\frac{g(x)}{h(x)}=1\ ,\hbox{ 
respectively }\hbox{limsup}_{x\rightarrow\alpha}\frac{g(x)}{h(x)}\leq 1\ ,$$
where $g(x),\ h(x)>0$.

{\bf Proof of Theorem \ref{majorationI}} 
We prove first a weaker statement assuming the
stronger inequality
$$2^{n-1}\tilde\Delta_{n-1}\sum_{k=1}^A
\left(1-k^2\left(2^{n-1}\tilde \Delta_{n-1}\frac{V_n}
{V_{n-1}}\sigma_n\right)^2\right)^{(n-1)/2}<1$$
where 
$$A=\lfloor \frac{2^{1-n} V_{n-1}}{V_n\tilde \Delta_{n-1}\sigma_n}
\rfloor.$$

Details for dealing with the extra factor 
$\left(\sum_{l\vert k}\frac{\mu(l)}{l^{n-1}}\right)$ will be 
given later.

Choose a positive real
number $\tilde\sigma_n<\sigma_n$ such that we have the equalities
$$A=\lfloor \frac{2^{1-n} V_{n-1}}{V_n\tilde \Delta_{n-1}\sigma_n}
\rfloor=\lfloor \frac{2^{1-n} V_{n-1}}{V_n\tilde \Delta_{n-1}
\tilde \sigma_n}
\rfloor$$
(where $\lfloor x\rfloor\in{\mathbf Z}$ denotes the integer part 
of $x\in{\mathbf R}$) and the inequality
$$2^{n-1}\tilde\Delta_{n-1}\sum_{k=1}^A\left(1-k^2
\left(2^{n-1}\tilde \Delta_{n-1}\frac{V_n}
{V_{n-1}}\tilde \sigma_n\right)^2\right)^{(n-1)/2}<
1.$$
We fix $\tilde \sigma_n$ in the sequel and introduce
$\epsilon=\frac{\sigma_n}{\tilde\sigma_n}-1>0$. 
We prove Theorem \ref{majorationI} for all $\mu$ huge enough
by showing the existence of a $\mu-$sequence $(s_0,\dots,s_{n-1},s_n)$
with $s_n\in I\cap {\mathbf N}$ where
$$I=[\tilde
\sigma_n\mu^{n/2}V_n,(1+\epsilon)\tilde\sigma_n\mu^{n/2}V_n]
=[\tilde
\sigma_n\mu^{n/2}V_n,\sigma_n\mu^{n/2}V_n]\ .$$ 
Since our computations rely on
strict inequalities involving volume-compu\-tations 
which are continuous in $\tilde
\Delta_{n-1}$, this will imply the weakened form (without the factor
$\left(\sum_{l\vert k}\frac{\mu(l)}{l^{n-1}}\right)$)
of Theorem \ref{majorationI}.

For $k=1,2,\dots\in{\mathbf N}$ we define finite subsets
$$I_k=\{s\in I\cap {\mathbf N} \vert\ \sum_{i=0}^{n-1}s_ix_i=ks 
\hbox{ for some }(x_0,\dots,x_{n-1})
\in B_{<\sqrt{\mu-k^2}}\cap {\mathbf Z}^n\}$$
of natural integers in $I\cap {\mathbf N}$ where
$B_{<\sqrt{\mu-k^2}}\cap {\mathbf Z}^n$ denotes the set of 
all integral vectors $(x_0,\dots,x_{n-1})\in {\mathbf Z}^n$ having
(squared Euclidean) norm strictly smaller than $\mu-k^2$.

We have
$$\begin{array}{c}
\displaystyle \vert\sum_{i=0}^{n-1}s_ix_i\vert\leq 
\sqrt{\sum_{i=0}^{n-1}s_i^2}\sqrt{\sum_{i=0}^{n-1}x_i^2}
\leq_{\mu\rightarrow\infty}\sqrt{\frac{\mu^{n-1}V_{n-1}^2}{4^{n-1}\tilde
    \Delta_{n-1}^2}}\sqrt{\mu-k^2}\\
\displaystyle \sim_{\mu\rightarrow\infty}
\frac{2^{1-n}\ \mu^{n/2}V_{n-1}}{\tilde\Delta_{n-1}}\end{array}$$
for $(x_0,\dots,x_{n-1})\in B_{<\sqrt{\mu-k^2}}$. 
This shows $I_k=\{\emptyset\}$ 
if 
$$k\geq A+1> \frac{2^{1-n}V_{n-1}}{V_n\tilde\Delta_{n-1}\tilde\sigma_n}
\ .$$

An extension $(s_0,\dots,s_{n-1},s_n)$ with $s_n\in I$ 
of a $\mu-$sequence $(s_0,\dots,s_{n-1})$ is a $\mu-$sequence if 
and only if $s_n\not\in \bigcup_{k=1}^A I_k$.

Introducing the sets 
$$X_k(a)=\{(x_0,\dots,x_{n-1})\in{\mathbf Z}^n\ \vert\ 
\sum_{i=0}^{n-1}s_ix_i\in kI\cap (k{\mathbf N}+a),
\ \sum_{i=0}^{n-1} x_i^2<\mu-k^2\},$$
we have obviously $\sharp(I_k)\leq \sharp(X_k(0))$.
This ensures the existence of a $\mu-$sequence 
$(s_0,\dots,s_{n-1},s_n)$
with $s_n\in I
\cap {\mathbf N}$ if we have
\begin{equation}
\sum_{k=1}^A \sharp(X_k(0)) 
<\sharp\{I\cap {\mathbf
  N}\}.
\label{fundineq}
\end{equation}
 
Denoting by 
$$X_k(*)=\{(x_0,\dots,x_{n-1})\in{\mathbf Z}^n\ \vert\ 
\frac{1}{k}\sum_{i=0}^{n-1}s_ix_i\in I,\ \sum_{i=0}^{n-1}
x_i^2<\mu-k^2\}$$
the union of the disjoint sets $X_k(0),X_k(1),\dots,X_k(k-1)$,
the following asymptotic equalities hold.

\begin{lem}\label{lemmaXk} We have 
$$\sharp(X_k(j))\sim_{\mu\rightarrow\infty} \frac{1}{k}\sharp(X_k(*))$$
for $j=0,\dots,k-1$. 
\end{lem}

It is thus enough to compute $\sharp(X_k(*))$ in order to get
an asymptotic estimation of $X_k(0)\sim_{\mu\rightarrow\infty}
\frac{1}{k}\sharp(X_k(*))$. We have
$$\begin{array}{l}
\displaystyle\sharp(X_k(*))=
\sharp\{(x_0,\dots,x_{n-1})\in{\mathbf Z}^n\ \vert\ 
\frac{1}{k}\sum_{i=0}^{n-1}s_ix_i\in I,\ \sum_{i=0}^{n-1} 
x_i^2<\mu-k^2\}\\
\quad \displaystyle\sim_{\mu\rightarrow\infty}
\hbox{Vol}\{(t_0,\dots,t_{n-1})\in {\mathbf E}^n\ \vert \ 
\sum_{i=0}^{n-1} t_i^2\leq \mu,\ \frac{1}{k}\sum_{i=0}^{n-1}s_it_i\in I
\}
\end{array}$$
and the requirement $\frac{1}{k}\sum_{i=0}^{n-1}s_it_i\in I$
amounts to the inequalities
$$k\tilde\sigma_n\mu^{n/2}V_n\leq \sum 
s_it_i\leq k\sigma_n\mu^{n/2}V_n.$$

For huge $\mu$ (and fixed $k$), 
the number $k\sharp(X_k)$ is thus essentially the volume 
$W_k$ of a subset of the $n-$dimensional ball of radius $\sqrt\mu$.
More precisely, this subset is delimited 
by the two parallel affine hyperplanes orthogonal to
$(s_0,\dots,s_{n-1})$ which are at distance
$$k\tilde D=k \frac{\tilde \sigma_n\mu^{n/2}V_n}
{\sqrt{\sum_{i=0}^n s_i^2}}\sim_{\mu\rightarrow\infty}
k\sqrt\mu\ 2^{n-1}\frac{V_n}{V_{n-1}}\tilde\Delta_{n-1}\tilde
\sigma_n$$
and $(1+\epsilon)k\tilde D$ of the origin.
 
We have thus
$$\begin{array}{lcl}
\displaystyle W_k&=&
\int_{k\tilde D}^{k(1+\epsilon)\tilde D}\left(\mu-t^2\right)^{(n-1)/2}dtV_{n-1}
\leq\epsilon k\tilde D
\left(\mu-k^2\tilde D^2\right)^{(n-1)/2}V_{n-1}\\
\displaystyle &\leq_{\mu\rightarrow\infty}&
\epsilon k\tilde \sigma_n\mu^{n/2}\ 2^{n-1}V_n\tilde
\Delta_{n-1}\left(1-k^2\left(2^{n-1}\tilde\sigma_n\frac{V_n}{V_{n-1}}
\tilde \Delta_{n-1}\right)^2\right)^{(n-1)/2}.\end{array}
$$
Using the asymptotic equalities
$\sharp(X_k)\sim_{\mu\rightarrow\infty}
\frac{W_k}{k}$, we get
$$\sum_{k=1}^A\sharp(X_k)\leq_{\mu\rightarrow\infty} 
\epsilon\tilde\sigma_n\mu^{n/2}\ 2^{n-1}\tilde\Delta_{n-1}V_n
\sum_{k=1}^A
\left(1-k^2\left(2^{n-1}\tilde \sigma_n\tilde\Delta_{n-1} \frac{V_n}
{V_{n-1}}\right)^2\right)^{(n-1)/2}.$$
Together with the obvious estimation
$$\sharp\{I\cap {\mathbf N}\}\sim_{\mu\rightarrow\infty}
\epsilon\tilde \sigma_n\mu^{n/2}V_n,$$
we have now
$$\frac{\sharp\{I\cap{\mathbf N}\}}{\sum_{k=1}^A\sharp(I_k)}
\geq_{\mu\rightarrow\infty}
\frac{2^{1-n}}{\tilde \Delta_{n-1}\sum_{k=1}^A
\left(1-k^2\left(2^{n-1}\tilde \sigma_n\tilde\Delta_{n-1} \frac{V_n}
{V_{n-1}}\right)^2\right)^{(n-1)/2}}>1$$
by assumption on the choice of $\tilde\sigma_n$. 
This proves the weak version (without the factor $\sum_{l\vert k}
\frac{\mu(l)}{l^{n-1}}$) of Theorem \ref{majorationI}
by inequation (\ref{fundineq})
since $\sharp\{I\cap{\mathbf N}\}\longrightarrow\infty$
if $\mu\rightarrow\infty$.

We consider now intersections
among the sets $I_1,I_2,\dots,I_A$ in order to deal 
with the factor $\sum_{l\vert k}
\frac{\mu(l)}{l^{n-1}}$ . This leads to a slightly better
estimation of $\sharp(\bigcup_{k=1}^A I_k)$ and completes the proof
of Theorem \ref{majorationI}.
 
Call an element $x=(x_0,\dots,x_{n-1})\in X_k(0)$ primitive if
it is not of the form $h{\mathbf Z}^n$ for an integral divisor
$h>1$ of $k$.
Call $x$ imprimitive otherwise. An imprimitive element
is of the form $h\tilde x$ with $\tilde x\in X_{k/h}(0)$ 
and contributes a common integer to the sets $I_k$ and $I_{k/h}$.
This implies the inequality
$$\sharp(\bigcup_{k=1}^A I_k)\leq \sum_{k=1}^A\sharp(X_k(0)_p)$$
where $X_k(0)_p\subset X_k(0)$ denotes the set of all primitive
elements in $X_k(0)$.

It is thus enough to estimate the number of primitive 
elements in $X_k(0)$. We have
$$\sharp(X_k(*)\cap h{\mathbf Z}^n)\sim_{\mu\rightarrow\infty}
\frac{1}{h^n}\sharp(X_k(*)).$$
We have obviously $X_k(a)\cap h{\mathbf Z}^n=\emptyset$
for $a\not\in h{\mathbf Z}$. Applying
Lemma \ref{lemmaXk}, obviously modified, to the sublattice 
$h{\mathbf Z}^n\subset{\mathbf Z}^n$ of index $h^n$ shows
$$\sharp(X_k(\alpha h)\cap h{\mathbf Z}^n)
\sim_{\mu\rightarrow\infty} \frac{1}{k/h}
\sharp(X_k(*)\cap h{\mathbf Z}^n)$$
for $\alpha=0,1,\dots,\frac{k}{h}-1$.
We get thus 
$$\sharp(X_k(0)\cap h{\mathbf Z}^n)\sim_{\mu\rightarrow\infty}
\frac{1}{kh^{n-1}}\sharp(X_k(*))\sim_{\mu\rightarrow\infty}
\frac{1}{h^{n-1}}\sharp(X_k(0)).$$
Since an element $x\in X_k(0)\cap h{\mathbf Z}^n$ belongs also
to $X_k(0)\cap l{\mathbf Z}^n$ for any natural divisor $l$ of $h$
and since $\sum_{l\vert h}\mu(l)=0$ for $h\geq 2$, the number
$\sharp(X_k(0)_p)$ 
of primitive elements in $X_k(0)$ is asymptotically given by
$$\left(\sum_{l\vert k}\frac{\mu(l)}{l^{n-1}}\right)\sharp(X_k(0)).$$
This leads to the majoration
$$\sharp(\bigcup_{k=0}^A I_k)\leq_{\mu\rightarrow\infty}\sum_{k=1}^A
\left(\sum_{l\vert k}\frac{\mu(l)}{l^{n-1}}\right)\sharp(X_k(0))$$
and proves Theorem \ref{majorationI}.
\hfill$\Box$


{\bf Proof of Lemma \ref{lemmaXk}}
The statement of Lemma \ref{lemmaXk} is equivalent to the 
asymptotic equalities
$$\frac{\sharp(X_k(j))-\sharp(X_k(i))}{
\sharp(X_k(*))}\sim_{\mu\rightarrow\infty}0$$
for $0\leq i,j<k$.

Fix $0\leq i<j<k$. Associate to an element 
$(x_0,x_1,\dots,x_{n-1})
\in X_k(j)$ the element $(x_0+i-j,x_1,\dots,x_{n-1})$
provided that it belongs to $X_k(i)$. This induces a bijection 
between subsets $\tilde X_k(j)$ and $\tilde X_k(i)$
of $X_k(j),X_k(i)$. The set of ``bad'' points 
$$B_k(i,j)=\left(X_k(j)\setminus \tilde X_k(j)\right)\cup 
\left(X_k(i)\setminus \tilde X_k(i)\right)$$
consists of some integral points at bounded 
Euclidean distance $<k\leq A$ from 
the boundary $\partial Z_k$ of the the set
$$Z_k=\{(z_0,\dots,z_{n-1})\in {\mathbf R}^n\  \vert 
\sum_{i=0}^{n-1}s_iz_i\in kI,\ \sum_{i=0}^{n-1} z_i^2\leq \mu-k^2\}_ .$$
This shows that 
$$\vert\sharp(X_k(j))-\sharp(X_k(i))\vert\leq
\sharp(B_k(i,j))\leq \hbox{vol}\left(N_{k+\sqrt{n}/2}(\partial Z_k)\right)
\sim O(\mu^{n-1})$$
where $N_{k+\sqrt{n}/2}(\partial Z_k)\subset {\mathbf R}^n$ denotes 
the set of all points at distance $\leq k+\sqrt{n}/2$ from the 
boundary $\partial Z_k$ of $Z_k$.

Since $\sharp(X_k(*))=O(\mu^n)$ this proves Lemma \ref{lemmaXk}.
\hfill $\Box$

\section{Proof of Theorem \ref{Minkowskibound}}
\label{secMink}

Although an intuitively correct proof using Theorem 
\ref{Hermiteineq} of Theorem \ref{Minkowskibound} is
easy, a rigorous proof is somewhat tedious and is the content of this 
section.

\subsection{Two auxiliary results}

The main ingredient for proving Theorem  \ref{Minkowskibound} is
the following result which says essentially that an attracting
fixpoint of a dynamical system is structurally stable.

\begin{prop} \label{fnlimit} Given a real interval $A\subset 
{\mathbf  R}$,
let $f_1,f_2,\dots:A\longrightarrow A$
be a sequence of uniformly converging functions with 
continuous and differentiable limit $f(x)=
\hbox{lim}_{n\rightarrow \infty}f_n(x)$ on $A$.
Suppose that $f$ has a fixpoint $\xi=f(\xi)\in A$ and 
suppose that we have $\hbox{sup}_{x\in A}\vert f'(x)\vert=\lambda<1$.

Then the sequence $s_n(x)$ of functions defined recursively 
by $s_0(x)=x$ and $s_n(x)=f_n(s_{n-1}(x))$ converges pointwise to the
constant function $\xi$.
\end{prop}

{\bf Proof} Given $\delta>0$ there exists an integer 
$N$ such that $\vert f_n(x)-f(x)\vert<\delta(1-\lambda)$ for all 
$x\in A$ and for all $n>N$. We have then for $m>N$
$$\begin{array}{lcl}
\displaystyle \vert s_m(x)-\xi\vert&=&
\displaystyle \vert f_m(s_{m-1}(x))-\xi\vert\\
&<&\delta(1-\lambda)+\vert f(s_{m-1}(x))-\xi\vert\\
&<&\delta-\lambda\delta+\lambda\vert s_{m-1}(x)-\xi\vert\\
&<&\delta-\lambda^2\delta+\lambda^2\vert s_{m-2}(x)-\xi\vert\\
& &\quad \vdots\\
&<&\delta-\lambda^{m-N}\delta+\lambda^{m-N}\vert
s_N(x)-\xi\vert\ .\end{array}$$
This shows $\vert s_m(x)-\xi\vert <2\delta$ if $$m>\hbox{max}(N,
N+\hbox{log}\left(\frac{\delta}{\vert s_N(x)-\xi\vert}\right)/\hbox{log}(\lambda))$$
and implies the result since we can choose $\delta>0$ 
arbitrarily small.\hfill $\Box$

\begin{rem} \label{fnremark} (i) The proof of Proposition 
\ref{fnlimit} shows in fact
$$\vert s_n(x)-\xi\vert\leq \lambda^n\vert x-\xi\vert +\sum_{k=1}^n
\lambda^{n-k}\hbox{sup}_{x\in A}\vert f_k(x)-f(x)\vert\ .$$

Asymptotically, we have thus $\vert s_n(\xi)-\xi\vert =O(
\hbox{sup}_{x\in N(\xi)}\vert f_n(x)-f(x)\vert)$ (where $N(\xi)\subset
A$ is an arbitrarily small fixed neighbourhood of $\xi$) if  
$\hbox{sup}_{x\in N(\xi)}\vert f_n(x)-f(x)\vert$  is
decreasing at a slower rate than powers of $\lambda$. 

\ \ (ii) If the sequence $f_n(x)=F(x,1/n)$ satisfies
the hypotheses of Proposition \ref{fnlimit} with $F(x,y)$ having
continuous partial derivatives of all orders up to $k+1$ 
in a neighbourhood of $(\xi,0)$, then there exist constants 
$a_1,a_2,\dots,a_k$ such that
$$s_n(x)=\xi+\frac{1}{1-\partial/\partial_x F}\left(\sum_{j=1}^k
  \frac{a_j}{j!} n^{-j}\right)+
O(n^{-(k+1)})$$
where $\frac{\partial^{a+b}F}{\partial x^a\partial y^b}$
denotes the obvious partial derivative of $F(x,y)$, 
evaluated at $(\xi,0)$. The formulae for the 
first three coefficients $a_1,a_2,a_3$
are
$$\begin{array}{lcl}
a_1&=&\displaystyle \partial/\partial_y F\\
a_2&=&\displaystyle 2a_1+\left(a_1\partial/\partial_x+
\partial/\partial_y\right)^2F\\
a_3&=&\displaystyle 12a_2-6a_1+6\partial/\partial_y\left(a_1\partial/\partial_x+\partial/\partial_y\right)F\\
&&\displaystyle \qquad +\left(a_1\partial/\partial_x+\partial/\partial_y\right)^3F\end{array}$$
In particular, for $F(x,y)$ analytic and non-constant in $y$, the 
sequence $s_n(x)$ is asymptotically independent from $x$.
\end{rem}

For $x\in(0,\infty)$ we consider the real analytic positive
function 
$$\tau(x)=\sum_{k=1}^\infty e^{-\pi(k/x)^2}=
\frac{1}{2}\theta_3(\frac{i}{x^2})-\frac{1}{2}\ ,$$
related to the third Jacobi-theta function 
$\theta_3(z)=\sum_{k\in{\mathbf Z}}e^{i\pi k^2z}$, cf. 
for instance Equation (6), page 102 in \cite{CS}. 
For $x>0$, we have $\tau'(x)=\frac{2\pi}{x^3}\sum_{k=1}^\infty
k^2e^{-\pi(k/x)^2}>0$ and the easy inequalities
$$\frac{x}{2}-1<\int_0^\infty
e^{-\frac{\pi}{x^2}t^2}dt-\int_0^1e^{-\frac{\pi}{x^2}t^2}dt<
\sum_{k=1}^\infty e^{-\pi(k/x)^2}<\int_0^\infty
e^{-\frac{\pi}{x^2}t^2}dt=\frac{x}{2}$$
for $x>0$ imply that $x\longmapsto \tau(x)$ is an increasing 
analytic diffeomorphism of $(0,\infty)$. 
The equation 
$$\frac{1}{x}=\tau\left(\frac{\Omega(x)}{x}\right)=\sum_{k=1}^\infty
e^{-k^2\pi(x/\Omega(x))^2}$$
defines thus a real positive analytic function $\Omega:(0,\infty)
\longrightarrow{\mathbf R}$.
Equivalently, the function $\Omega$ is given by $\Omega(x)=x
\psi\left(\frac{1}{x}\right)$ where the analytic diffeomorphism 
$\psi$ satisfies $\psi(\tau(x))=\tau(\psi(x))=x$ for all $x>0$
and is the reciprocal function of $\tau$.

The proof of Theorem \ref{Minkowskibound}
uses the following result.

\begin{prop} \label{Omegamonotone} The application
$$x\longmapsto \Omega(x)=x\psi\left(\frac{1}{x}\right)$$
defines a continuous map from $(0,\infty)$ onto $(2,\infty)$
which is strictly increasing for $x>2$. 
It has a unique fixpoint 
$\xi=\frac{1}{\tau(1)}=\frac{1}{\sum_{k=0}^\infty e^{-\pi k^2}}
\sim 23.13882534$ which is attracting under iteration since
$$\Omega'(\xi)=1-\frac{\tau(1)}{\tau'(1)}
=1-\frac{\sum_{k=1}^\infty e^{-\pi k^2}}{2\pi\sum_{k=1}^\infty
  k^2e^{-\pi k^2}}\sim 0.9135652<1\ .$$
\end{prop}

\subsection{Proof of Theorem \ref{Minkowskibound}}

Given an $(n-1)-$dimensional lattice of density $\tilde \Delta_{n-1}$,
Theorem \ref{majorationI} implies the existence of an $n-$dimensional
lattice with density $\tilde \Delta_n$ arbitrarily close to 
$\frac{1}{2^n\overline \sigma}$ for $\overline \sigma>0$ defined
by
$$2^{n-1}\tilde \Delta_{n-1}\sum_{k=1}^A
\sqrt{1-k^2\left(2^{n-1}\tilde \Delta_{n-1}\frac{V_n}{V_{n-1}}
\overline{\sigma}\right)^2\quad}^{n-1}=1$$
where
$$A=\left\lfloor \frac{2^{1-n}V_{n-1}}{V_n\tilde\Delta_{n-1}
  \overline{\sigma}}\right\rfloor\ .$$ 
Given a positive constant $\epsilon>0$ and a natural integer
$N\in{\mathbf N}$,
there exists thus a sequence of lattices $\Lambda_1,\Lambda_2,\dots,
\Lambda_N$ of dimensions $1,2,\dots,N$ with densities $\tilde \Delta_1=1,
\tilde\Delta_2,\dots,\tilde\Delta_N$ satisfying
$$\tilde\Delta_m\geq (1-\epsilon)\frac{d_m}{2^m},\ m=1,\dots,N$$
where $d_1=2$ and $d_2,d_3,\dots,d_N$ are recursively defined by
the equation 
$$d_{n-1}\sum_{k=1}^{A_n}
\sqrt{1-k^2\left(\frac{d_{n-1}}{d_n}\ \frac{V_n}{V_{n-1}}\right)^2
\quad }^{n-1}=1\hbox{ with }
A_n=\left\lfloor \frac{d_n\ V_{n-1}}{d_{n-1}\ V_n}\right\rfloor\ .$$
Equivalently, the sequence $d_1,d_2,\dots$ is given by 
$d_1=2,d_2=f_1(2),d_3=f_2(d_2),\dots,d_{n+1}=f_n(d_n),\dots$
where $f_1,f_2,\dots:(0,\infty)\longrightarrow(0,\infty)$
are the functions defined implicitely by the equations
$$x\sum_{k=0}^{\lfloor f_n(x)V_n/(xV_{n+1})\rfloor}
\sqrt{1-k^2\left(\frac{xV_{n+1}}{f_n(x)V_n}\right)^2\quad }^n=1\ .$$

Stirlings formula $n!=\sqrt{2\pi
  n}(n/e)^n(1+O(1/n))$ shows 
$$V_{n+1}/V_{n}=\sqrt{\pi}\frac{(n/2)!}{((n+1)/2)!}
=\sqrt{2\pi/n}\ \left(1+O(1/n)\right)\ .$$ 
We have thus asymptotically
$$\begin{array}{cl}
\displaystyle 
1=&\displaystyle x\sum_{k=0}^{\lfloor f_n(x)V_n/(xV_{n-1})\rfloor}
\sqrt{1-k^2\left(\frac{xV_{n+1}}{f_n(x)V_n}\right)^2\quad }^n\\
\displaystyle\phantom{1}=&\displaystyle
\left(x\sum_{k=1}^\infty e^{-k^2\pi\left(x/f_n(x)\right)^2}\right)(1+
O(1/n))\end{array}$$
and $f_n(x)\longrightarrow \Omega(x)$ uniformly on any
compact subset of $(0,\infty)$. By Proposition \ref{Omegamonotone}
we can find $\alpha<\xi=\left(\sum_{k=1}^\infty
  e^{-\pi k^2}\right)^{-1}\sim 23.14<\beta$ such that $\Omega'(x)\leq 19/20$
for $x\in [\alpha,\beta]$.
We have thus uniform convergency $f_n(x)\longrightarrow \Omega(x)$
for $x\in [\alpha,\beta]$, and there exists an integer $N_\xi$
such that $f_n([\alpha,\beta])\subset
[\alpha,\beta]$ for all $n\geq N_\xi$. 
Proposition \ref{fnlimit} shows now 
$$\hbox{lim}_{n\rightarrow\infty}d_n=\xi$$
which ends the proof. \hfill $\Box$

The following Table illustrates the convergence of the 
sequence $d_1=2,d_2=f_1(d_1),\dots$:

$$\begin{array}{rrrr}
    1& 2.00000000& 2.00000000&  0\\
    2& 3.62759873& 3.99997210& -0.7447467\\
    4& 8.08369319& 7.92472241&  0.6358831\\
    8&18.71971890&14.38756801& 34.6572071\\
   16&30.69030131&20.71395996&159.6214617\\
   32&29.45114255&22.98242063&206.9991014\\
   64&25.53248635&23.13821340&153.2334688\\
  128&24.17810739&23.13882533&133.0281029\\
  256&23.63011883&23.13882534&125.7711333\\
  512&23.37820694&23.13882534&122.5633803\\
 1024&23.25703467&23.13882534&121.0463495\\
\end{array}$$

The first column shows the indices $n$, choosen as successive 
powers of $2$. The second column shows the corresponding value
of $d_n$. The third column shows the $(n-1)-$th iteration of $\Omega$,
starting from the initial value $2$. The last column 
is the difference between the second and third column, multiplied
by $n$ and illustrates the expected finer asymptotic properties.

Asymptotically, the number $d_n$ is roughly given by
$$23.13882534+119.58193\frac{1}{n}+1473.8282\frac{1}{n^2}+
25774.448\frac{1}{n^3}+\dots$$
(cf assertion (ii) of Remark \ref{fnremark}). 

\subsection{Proof of Proposition \ref{Omegamonotone}}

Using the orientation-reversing diffeomorphism $x=\frac{1}{\tau(Y)}
\longmapsto
Y=\psi\left(\frac{1}{
x}\right)$ of $(0,\infty)$ we have
$$\frac{Y}{\tau(Y)}=\Omega(x)=\Omega\left(\frac{1}{\tau(Y)}\right)
\ .$$

The inequality $\tau(Y)<\frac{Y}{2}$ shows $\Omega(x)=\frac{Y}
{\tau(Y)}>2$ and $\frac{2Y}{Y-2}>\frac{Y}{\tau(Y)}$
implies $\hbox{lim}_{x\rightarrow 0_+}\Omega(x)=2$.
Since 
$$\hbox{lim}_{Y\rightarrow 0_+}\frac{Y}{\tau(Y)}=
\hbox{lim}_{Y\rightarrow 0_+}Ye^{\pi^2/Y^2}
\left(1+\sum_{k=2}^\infty e^{-\pi(k^2-1)/Y^2}\right)^{-1}=\infty$$
the map $\Omega$ is a surjection onto $(2,\infty)$.

Since $x\longmapsto Y$ is orientation reversing,
$\frac{d}{dx}\Omega(x)>0$ is equivalent to strict positivity of 
$$Y^2\frac{d}{dY}\left(\frac{\tau(Y)}{Y}\right)=Y\tau'(Y)-\tau(Y)=
\frac{1}{Y^2}\sum_{k=1}^\infty\left(2\pi k^2-Y^2\right)
e^{-\pi(k/Y)^2}$$
which obviously holds for $Y\leq \sqrt{2\pi}$ corresponding
to 
$$x\geq \frac{1}{\tau(\sqrt{2\pi})}=\frac{1}{\sum_{k=1}^\infty
e^{-\pi k^2/2}}\sim 1.38\ .$$
This implies that $\Omega$ restricts to an increasing  diffeomorphism
from $(2,\infty)$ onto $(\Omega(2),\infty)$ and since
$\Omega(x)>2$, the map $x\longmapsto \Omega(x)$
has a unique fixpoint at $\xi=\frac{1}{\tau(1)}$.
\hfill $\Box$

\section{Final remarks}\label{sectioncomments}

The inequality
$$\sharp(\bigcup_{k=1}^A I_k)\leq\sum_{k=1}^A\sharp(X_k(0)_p)$$
appearing in the proof of Theorem \ref{majorationI}
is probably not sharp. A smaller upper bound for the cardinality 
$\sharp(\bigcup_{k=1}^A I_k)$ would thus improve the results of 
this paper.

The inequality above can be decomposed into the two inequalities
$$\sharp(\bigcup_{k=1}^AI_k)\leq\sum_{k=1}^A\sharp(I_{k,p})$$
and
$$\sharp(I_{k,p})\leq \sharp(X_k(0)_p)$$
where we denote by $I_{k,p}\subset I_k$ the subset of integers
corresponding to primitive elements. If the subsets 
$I_{1,p},\dots,I_{A,p}$ are asymptotically
``independent'' in the sense that 
$$\sharp(\bigcap_{j=1}^l
I_{k_j,p})/\sharp(I\cap{\mathbf Z})\sim_{\mu\rightarrow\infty}
\prod_{j=1}^l\big(\sharp(I_{k_j,p})/\sharp(I\cap{\mathbf Z})\big)\ ,$$
for $\{I_{k_1,p},\dots,I_{k_l,p}\}\subset\{I_{1,p},\dots,I_{A,p}\}$ 
a subset of $l$ distinct elements,
one can neglect the contributions corresponding to
$k=2,\dots,A$. This would lead to a small improvement.

A probably much more important improvement would result from a better
understanding of the inequality $\sharp(I_{k,p})\leq
\sharp(X_k(0)_p)$.

Instead of working with sublattices of ${\mathbf Z}^{n+1}$
orthogonal to a given vector $(s_0,\dots,s_n)\in {\mathbf Z}^{n+1}$, it
is possible to consider sublattices ${\mathbf Z}^{n+a}$
which are orthogonal to a set of $a\geq 2$ linearly independent
vectors in ${\mathbf Z}^{n+a}$. One might also replace the 
standard lattice ${\mathbf Z}^{n+1}$ by other lattices, e.g.
sublattices of dimension $n$ in ${\mathbf Z}^{n+1}$ (which approximate
homothetically an arbitrary lattice by Proposition \ref{Lestdense})
or of finite index in ${\mathbf Z}^{n+1}$.

Extending
finite $\mu-$sequences in an optimal way into longer $\mu-$sequences
amounts geometrically to the familiar process  
of lamination for lattices (see for instance \cite{CS} or \cite{M}). 
The existence of an integer $s\in I\setminus I_1$
implies indeed the existence of a 
point $P\in {\mathbf E}^{n-1}$ which is
far away from any lattice point of the affine lattice 
$\{(x_0,\dots,x_{n-1})\ \vert\sum x_is_i=s\}\subset{\mathbf Z}^n$ 
and corresponds thus to a
\lq\lq hole'' of the lattice.

The present version of this paper ows much to Fedor Petrov 
whose pertinent questions clarified and improved (and changed the
title of) a preliminary version, see \cite{Ba}.

I thank also A. Marin, 
J. Martinet, P. Sarnak, B. Venkov and J-L. Verger-Gaugry
for helpful comments and interest in  this work.

Roland Bacher, INSTITUT FOURIER, Laboratoire de Math\'ematiques, 
UMR 5582 (UJF-CNRS), BP 74, 38402 St MARTIN  D'H\`ERES Cedex (France),
e-mail: Roland.Bacher@ujf-grenoble.fr

\end{document}